\documentclass[12pt]{amsart}

\usepackage{color}
\usepackage{amsbsy}
\usepackage{amssymb,amsmath,amstext,mathrsfs}
\usepackage{graphicx}

\newtheorem{theorem}{Theorem}[section]
\newtheorem{lemma}[theorem]{Lemma}
\newtheorem{proposition}[theorem]{Proposition}
\newtheorem{corollary}[theorem]{Corollary}
\theoremstyle{definition}
\newtheorem{definition}[theorem]{Definition}

\newtheorem{remark}[theorem]{Remark}
\newtheorem{claim}{Claim}

\newcommand{\va}{\left|}
\newcommand{\vb}{\right|}
\newcommand{\vd}{\right\| }
\newcommand{\vc}{\left\| }

\newcommand{\pl}{\left(}
\newcommand{\pr}{\right)}
\newcommand{\ql}{\left[}
\newcommand{\qr}{\right]}
\newcommand{\tl}{\left\{}
\newcommand{\tr}{\right\}}

\newcommand{\ra}{\rightarrow}

\begin{document}

\title[Biseparating maps]{\bf Biseparating maps
between Lipschitz function spaces}

\author{Jes\'us Araujo and Luis Dubarbie}

\address{Departamento de Matem\'aticas, Estad\'istica y
Computaci\'on, Facultad de Ciencias, Universidad de Cantabria,
Avenida de los Castros s/n, E-39071, Santander, Spain.}

\email{araujoj@unican.es} \email{dubarbiel@unican.es}

\thanks{Research partially supported by the Spanish Ministry of
Science and Education (MTM2006-14786). L. Dubarbie was partially supported by a predoctoral
grant from the University of Cantabria and the Government of Cantabria. }



\begin{abstract}
For complete metric spaces $X$ and $Y$,  a description of
linear biseparating maps between spaces of vector-valued Lipschitz
functions defined on $X$ and $Y$ is provided. In particular it is proved that
 $X$ and $Y$ are bi-Lipschitz homeomorphic, and the automatic continuity of
such maps is derived in some cases. Besides, these results are used to
characterize the separating bijections between scalar-valued Lipschitz
function spaces when $Y$ is compact.
\end{abstract}

\maketitle

\section{Introduction}

Separating  maps, also called disjointness preserving maps,
between spaces of scalar-valued continuous functions defined on compact or
locally compact spaces have drawn the attention of researchers in
last years (see for instance  \cite{FH}, \cite{J},   \cite{JW1} and \cite{KN}). Roughly speaking, a (bijective) linear operator $T$ between two spaces of functions
is said to be separating if $(Tf) \cdot (Tg) =0$ whenever $f \cdot g=0$ (see Definition~\ref{imse}).

 No results are known so far for the case when the map is defined between spaces of Lipschitz functions, even if successful attempts have been made for some special subalgebras. Namely,  
Jim\'enez-Vargas recently obtained  the representation  of separating maps  defined between {\em little}
Lipschitz algebras  on {\em compact} metric spaces (see \cite{JV}). Unfortunately proofs rely heavily on the properties of these algebras and on the compactness of spaces, so that they cannot carry over to the general case.
Also,  in the recent paper \cite{GJ2}, Garrido and Jaramillo study a related problem: find those metric spaces $X$ for which the algebra of bounded Lipschitz functions on $X$ determines the Lipschitz structure of $X$. But even if separating maps are related with algebra isomorphisms, their techniques cannot be used here either.

The aim of this paper is to study such maps and obtain their general representation. In fact, we do not restrict ourselves to the scalar setting and we deal with the vector-valued case as well. As usual,
when  spaces of functions taking values in arbitrary normed spaces are involved, the condition for an operator of being separating is not enough 
to ensure  a good representation, and we must require  the inverse map to be  separating too (see for instance
\cite{AK}, \cite{A1}, \cite{A2}, \cite{A3},
\cite{AJ}, \cite{GJW}, \cite{HBN}, \cite{JW}; see also \cite[Theorem 5.4]{A4} and \cite{DU} for special cases where this may not be true). We also drop any requirement of compactness on the 
metric spaces where functions are defined, and completeness is assumed instead.

Other papers where related operators have been recently studied in similar contexts are \cite{BA},  \cite{GJ} and \cite{MCX} (see also \cite{S} and \cite{LS}).

\medskip

The paper is organized as follows. In Section~\ref{sect2} we give some definitions and notation that 
we use throughout the paper. In Section~\ref{maine} we state the main results. In Section~\ref{sect3}  we give  some  properties of  spaces
of Lipschitz functions 
that we use later.  Section~\ref{sect4} is devoted to prove the main results concerning 
biseparating maps between spaces of vector-valued Lipschitz
functions.
In particular, apart from obtaining their general form, we show that the
underlying spaces are bi-Lipschitz homeomorphic and, when $E$ and $F$ are complete, we obtain the automatic continuity of some related maps.
Finally, in Section~\ref{sect5} we prove that every bijective separating map between spaces of scalar-valued Lipschitz
functions defined on compact metric spaces is indeed biseparating.

\section{Preliminaries and notation}\label{sect2}

Let $(X, d_1)$ and $(Y, d_2)$ be metric spaces. Recall that a map $f:X \rightarrow Y$ is said to be \emph{Lipschitz} if there
exists a constant $k\geq 0$ such that
\begin{displaymath}
d_2 ( f(x) , f(y) ) \leq k \ d_1(x,y)
\end{displaymath}
for each $x,y\in X$. The least such $k$ is called the
\emph{Lipschitz number} of $f$ and will be denoted by $L(f)$.
Equivalently, $L(f)$ can be defined as
\begin{displaymath}
L(f):=\mathrm{sup}\left\{\frac{ d_2 (f(x) , f(y) )}{d_1 (x,y)}:x,y\in
X, x\neq y\right\}.
\end{displaymath}

When $f$ is bijective and  both $f$ and $f^{-1}$ are Lipschitz, we will say that $f$ is {\em bi-Lipschitz}.

If $E$ is  a  
$\mathbb{K}$-normed space,
where $\mathbb{K}$ stands for the field of real or complex numbers, then  
$\mathrm{Lip}(X,E)$ will denote the space of all \emph{bounded} $E$-valued
Lipschitz functions defined on $X$. If $E=\mathbb{K}$, then we put 
$\mathrm{Lip}(X):=\mathrm{Lip}(X,E)$.

It is well known that $\mathrm{Lip}(X,E)$ is a normed space  endowed with the norm
\begin{displaymath}
\|f\|_{L}=\mathrm{max}\left\{\|f\|_{\infty},L(f)\right\}
\end{displaymath}
for each $f\in \mathrm{Lip}(X,E)$ (where $\vc  \cdot \vd_{\infty}$ denotes
the usual supremum norm), which is  complete when $E$ is a Banach space.

\medskip
 
From now on, unless otherwise stated, we will suppose that $X$ and $Y$ are
{\em bounded} complete  metric spaces (see  Remark~\ref{helena}). In general, we will use $d$  to 
denote the metric in both spaces. 

 For
$x_{0}\in X$ and $r>0$, $B(x_{0},r)$ will denote the open ball
$\{x\in X:d(x,x_{0})<r\}$.
Finally, if $A$ is a subset of a topological space $Z$,
$\mathrm{cl}_{Z} A$ stands for the closure of $A$ in $Z$.

We will suppose that $E$ and $F$ are
$\mathbb{K}$-normed spaces. Given a function $f$ defined on $X$ and taking values on $E$, we
define the \emph{cozero set} of $f$ as $\mathrm{coz}(f):=\{x\in
X:f(x)\neq 0\}$. 
 Also, for each $\mathbf{e}\in E$,
$\widehat{\mathbf{e}}:X\rightarrow E$ will be the constant function taking
the value $\mathbf{e}$.
On the other hand, if $(f_{n})$ is a sequence of functions, then
$\sum_{n=1}^{\infty}f_{n}$ denotes its (pointwise) sum.

Finally, we will denote by $L'(E, F)$  the set of linear and bijective 
maps from $E$ to $F$, and by $L(E, F)$ the subset of  all {\em continuous} operators of $L'(E,F)$ .

\medskip

We now give the definition of separating and biseparating maps in the context of Lipschitz function spaces.

\begin{definition}\label{imse}
A linear map $T:\mathrm{Lip}(X,E)\rightarrow \mathrm{Lip}(Y,F)$ is
said to be \emph{separating} if $\mathrm{coz}(Tf)\cap
\mathrm{coz}(Tg)=\emptyset$ whenever $f,g\in \mathrm{Lip}(X,E)$
satisfy $\mathrm{coz}(f)\cap \mathrm{coz}(g)=\emptyset$.
Moreover, $T$ is said to be \emph{biseparating} if it is bijective and
both $T$ and $T^{-1}$ are separating.
\end{definition}

Equivalently, a map $T:\mathrm{Lip}(X,E)\rightarrow
\mathrm{Lip}(Y,F)$ is \emph{separating} if it is linear and
$\|Tf(y)\|\|Tg(y)\|=0$ for all $y\in Y$, whenever $f,g\in
\mathrm{Lip}(X,E)$ satisfy $\|f(x)\|\|g(x)\|=0$ for all $x\in X$.

\section{Main results}\label{maine}

Our first result gives a general description of biseparating maps.

\begin{theorem}\label{imse2008}
Let $T:\mathrm{Lip}(X,E)\rightarrow \mathrm{Lip}(Y,F)$ be a
biseparating map. Then there exist a bi-Lipschitz homeomorphism $h: Y \ra X$ and a map $J : Y \ra L'(E, F)$ such that
$$Tf(y)=(J y) (f(h(y)))$$
for all $f\in \mathrm{Lip}(X,E)$ and $y\in Y$.
\end{theorem}

Due to the representation given above, we see that when $T$ is continuous, then $J y$ belongs to $L(E,F)$ for every $y \in Y$. In particular we also have that, for $y, y' \in Y$ and $\mathbf{e} \in E$, the map
$ \vc T \widehat{\mathbf{e}} (y) - T \widehat{\mathbf{e}} (y') \vd \le \vc T \vd \vc \mathbf{e} \vd d(y, y')$.
Consequently, the map $y \in Y \mapsto J y \in L(E,F)$ is continuous when  $L(E,F)$ is endowed with the usual norm.

Of course Theorem~\ref{imse2008} does not give an answer to whether or not a biseparating map is necessarily  continuous. In fact, automatic continuity cannot be derived in general. Nevertheless, in some cases an associated continuous operator can be defined. This is done in Theorem~\ref{segundano}. We first give a result concerning
continuity of maps $J y$.

Given a biseparating map $T:\mathrm{Lip}(X,E)\rightarrow \mathrm{Lip}(Y,F)$, 
we denote $$Y_{d}:=\{y\in Y: J y \ \mathrm{is \ discontinuous} \}.$$

\begin{proposition}\label{kurda}
Let $T:\mathrm{Lip}(X,E)\rightarrow \mathrm{Lip}(Y,F)$ be a
biseparating map. Then the set $\tl \vc J y \vd : y \in Y \setminus Y_d \tr$ is bounded. Moreover, $Y_d$ is finite and each point of $Y_d$ is isolated in $Y$. 
\end{proposition}

An immediate consequence is the following.

\begin{corollary}\label{july11ensinka}
Let $T:\mathrm{Lip}(X,E)\rightarrow \mathrm{Lip}(Y,F)$ be a
biseparating map. If $X$ is infinite, then $E$ and $F$ are isomorphic.
\end{corollary}

Another immediate consequence  of Proposition~\ref{kurda} and Theorem~\ref{imse2008} is that $Y\setminus Y_{d}$ is 
complete, and that the restriction of $h$ to this set is a homeomorphism onto $X\setminus h(Y_{d})$. This allows us to
introduce in a natural way a new biseparating map defined in a related domain.

\begin{theorem}\label{segundano}
Suppose that $E$ and $F$ are complete. Let $T:\mathrm{Lip}(X,E)\rightarrow \mathrm{Lip}(Y,F)$ be a
biseparating map, and let $J$ and $h$ be as in Theorem~\ref{imse2008}. Then
$T_d :\mathrm{Lip}(X\setminus h(Y_{d}),E)\rightarrow \mathrm{Lip}(Y\setminus Y_{d},F)$, defined as 
$$T_d f(y) := (J y)(f(h(y)))$$
for all $f\in \mathrm{Lip}(X\setminus h(Y_{d}),E)$ and $y\in Y \setminus Y_{d}$,
is biseparating and continuous.
\end{theorem}

In the case when $Y$ is compact and we deal with spaces of scalar-valued   functions, the assumption on $T$ of being just separating and bijective is enough to obtain both its automatic continuity  and the fact that it is biseparating.

\begin{theorem}\label{gerria}
Let $T:\mathrm{Lip}(X)\rightarrow \mathrm{Lip}(Y)$ be a bijective
and separating map. If $Y$ is compact, then $T$ is biseparating and continuous.
\end{theorem}

\begin{remark}\label{helena}
Recall that we are assuming that the metrics in $X$ and $Y$ are bounded. Nevertheless results can be translated
to the case of unbounded metric spaces. Let $d_1$ be an unbounded metric in $X$ such that $(X, d_1)$ is complete. Then $d_1' := \min \tl 2, d_1 \tr$ is a
 {\em bounded} complete metric in $X$ and the topology induced by both metrics is the same.
 Following the same ideas as in \cite[Proposition 1.7.1]{W}, we can also see that the identity map of the space  $\mathrm{Lip}(X,E)$ (with respect to $d_1$) onto itself (with respect to $d_1'$) is an  isometric isomorphism. 
It is easy to see now that if 
 $d_2$ is a (bounded or unbounded) complete metric in $Y$, then a map $f: (Y, d_2) \ra (X, d_1')$ is Lipschitz if and only if 
$f: (Y, d_2) \ra (X, d_1)$
is  what is called
{\em Lipschitz in the small}, that is, there exist $r,  k>0$ such that
  $ d_1 (f(y), f(y') ) \le k \ d_2 (y, y')$ 
whenever $d_2 (y, y') < r$.
\end{remark}

\section{Lipschitz function spaces}\label{sect3}

Notice that since every complete metric space $X$ is
completely regular, it admits a Stone-\v{C}ech compactification,
which will be denoted by $\beta X$. Recall that this implies that every continuous map
$f:X\rightarrow \mathbb{K}$ can be extended to a continuous map
$f^{\beta X}$ from $\beta X$ into $\mathbb{K}\cup \{\infty\}$. In particular,
given a continuous map $f:X\rightarrow E$, we will denote by
$\|f\|^{\beta X}$ the extension of
$\|.\|\circ f:X\rightarrow \mathbb{K}\cup \{\infty\}$ to $\beta X$.

Now, we suppose that $A(X)$ is a subring of the space of continuous
functions $C(X)$ which separates each point of $X$ from each point
of $\beta X$. We introduce in $\beta X$ the equivalence
relation
\begin{displaymath}
x\sim y \Leftrightarrow f^{\beta X}(x)=f^{\beta X}(y)
\end{displaymath}
for all $f\in A(X)$. In this way, we obtain the quotient space
$\gamma X:=\beta X/\sim$, which is a new compactification of $X$.
 Besides, each $f\in
A(X)$ is continuously extendable to a map
$f^{\gamma X}$ from $\gamma X$ into $\mathbb{K}\cup \{\infty\}$.
In this context,  $A(X)$ is said to be 
{\em strongly regular} if given $x_0 \in \gamma X$ and a nonempty closed subset $K$ of $\gamma X$ that does not contain $x_0$, there exists $f \in A(X)$ such that $f^{\gamma X} \equiv 1$ on a neighborhood of $x_0$ and $f^{\gamma X} (K) \equiv 0$.

Finally, assume that $A(X,E)\subset C(X,E)$ is an $A(X)$-module. We
will say that $A(X,E)$ is \emph{compatible} with $A(X)$ if, for
every $x\in X$, there exists $f\in A(X,E)$ with $f(x)\neq 0$, and
if, given any points $x,y\in \beta X$ such that $x\sim y$, we have
$\|f\|^{\beta X}(x)=\|f\|^{\beta X}(y)$ for every $f\in A(X,E)$. In
this case, it is easy to see that $\|.\|\circ f:X\rightarrow
\mathbb{K}\cup \{\infty\}$ can be continuously extended to
$\|f\|^{\gamma X}$ from $\gamma X$ into $\mathbb{K}\cup \{\infty\}$.

It is straightforward to check that, if $f\in
\mathrm{Lip}(X)$ and $g\in \mathrm{Lip}(X,E)$, then $f\cdot g\in
\mathrm{Lip}(X,E)$, that is,
\begin{lemma}\label{trestres}
$\mathrm{Lip}(X,E)$ is a $\mathrm{Lip}(X)$-module.
\end{lemma}

\begin{remark}\label{macaya}
We introduce two families of Lipschitz functions that will be used 
later. Given $x_{0}\in X$ and $r>0$, the function
$\psi_{x_{0},r}:X\rightarrow \mathbb{K}$ defined as
\begin{displaymath}
\psi_{x_{0},r}(x):=\mathrm{max}\left\{0,1-\frac{d(x,x_{0})}{r}\right\}
\end{displaymath}
for all $x\in X$, belongs to $\mathrm{Lip}(X)$ and satisfies
$\psi_{x_{0},r}(x_{0})=1$,
$\mathrm{coz}(\psi_{x_{0},r})=B(x_{0},r)$,
$\|\psi_{x_{0},r}\|_{\infty}=1$, and $L(\psi_{x_{0},r})=1/r$.
On the other hand, another  Lipschitz function we will use  is
\begin{displaymath}
\varphi_{x_{0},r}(x):=\mathrm{max}\left\{0,1-\frac{d(x,B(x_{0},r))}{r}\right\}
\end{displaymath}
for all $x\in X$, which satisfies
$\varphi_{x_{0},r}(B(x_{0},r))\equiv 1$,
$\mathrm{coz}(\varphi_{x_{0},r})=B(x_{0},2r)$,
$\|\varphi_{x_{0},r}\|_{\infty}=1$, and $L(\varphi_{x_{0},r})=1/r$.
\end{remark}

Clearly, given $f\in \mathrm{Lip}(X,E)$, $\|.\|\circ f\in
\mathrm{Lip}(X)$. Then, by the definition of the 
equivalence relation $\sim$ in $\beta X$ given above and  the function $\psi_{x_{0},r}\in \mathrm{Lip}(X)$  for each
$x_{0}\in X$ (see Remark~\ref{macaya}), we obtain the next lemma.

\begin{lemma}\label{tres4}
$\mathrm{Lip}(X,E)$ is compatible with $\mathrm{Lip}(X)$.
\end{lemma}

\begin{lemma}\label{tresc}
$\mathrm{Lip}(X)$ is strongly regular.
\end{lemma}

\begin{proof}
Let $K$ and $L$ be two disjoint closed subsets of $\gamma X$. Since
$\gamma X$ is compact, there exists $f_{0}\in C(\gamma X)$, $0\leq
f_{0} \leq 1$, satisfying $f_{0}(K)\equiv 0$ and $f_{0}(L)\equiv 1$.
Obviously $K_0:=\{x\in \gamma X: f_{0}(x) \le 1/3\}$ and $L_0:=\{x\in \gamma
X: f_{0}(x) \ge 2/3\}$ are disjoint compact neighborhoods of $K$ and $L$,
respectively. Consider now $K_1 := K_0 \cap X$ and  $L_1 := L_0 \cap X$. We claim that $d(K_1,L_1)>0$.

Suppose this is not true, so for each $n\in \mathbb{N}$ there exist
$x_{n}\in K_1$ and $z_{n}\in L_1$ such that
$d(x_{n},z_{n})<1/n$. Since $K_0$ is compact,
$\{x_{n}:n\in \mathbb{N}\}$ has a limit point $x_{0}$ in $K_0$. 
Consequently, there exists a net $(x_{\alpha})_{\alpha\in \Omega}$
in $\{x_{n}:n\in \mathbb{N}\}$ which converges to $x_{0}$. Clearly,
for each $\alpha\in \Omega$, $x_{\alpha}=x_{n_{\alpha}}$ for some
$n_{\alpha}\in \mathbb{N}$. Next, we consider the net
$(z_{\alpha})_{\alpha\in \Omega}$ in $\{z_{n}:n\in \mathbb{N}\}$
defined, for each $\alpha\in \Omega$, by
$z_{\alpha}:=z_{n_{\alpha}}$ whenever $x_{\alpha}=x_{n_{\alpha}}$.
By the compactness of $L_0$, we know that there exists a subnet
$(z_{\lambda})_{\lambda\in \Lambda}$ of $(z_{\alpha})_{\alpha\in
\Omega}$ converging to a point $z_{0}$ in $L_0$. 

We are going
to prove that $x_{0}=z_{0}$, which is absurd because $K_0 \cap L_0 = \emptyset$. Obviously if $x_0$ or $z_0$ belongs to $X$, then we would have $x_0 = z_0$, so we assume that this is not the case. 
Let $U$ and $V$ be open neighborhoods of $x_0$ and $z_0$, respectively, and let $n_0 \in \mathbb{N}$. We are going to
 see that there exists $n \ge n_0$, $n \in \mathbb{N}$, such that $x_n \in U$ and $z_n \in V$. Without loss of generality we assume that $x_1 , \ldots, x_{n_0} \notin U$ and $z_1 , \ldots, z_{n_0} \notin V$. Since 
$ (x_{\lambda})_{\lambda\in \Lambda}$ and $(z_{\lambda})_{\lambda\in
\Lambda}$ converge to $x_{0}$ and $z_{0}$, respectively, there exist
$\lambda_{1}^{x_{0}}\in \Lambda$ and $\lambda_{1}^{z_{0}}\in
\Lambda$ such that $x_{\lambda}\in U$ for all
$\lambda\geq\lambda_{1}^{x_{0}}$ and $z_{\lambda}\in V$
 for all
$\lambda\geq\lambda_{1}^{z_{0}}$. Taking
$\lambda \in \Lambda$ such that $\lambda \geq \lambda_{1}^{x_{0}}, \lambda_{1}^{z_{0}}$, it is clear that $x_{\lambda}\in U$ and
$z_{\lambda}\in V$.
Now,  there exists $n_{\lambda}\in
\mathbb{N}$ such that $x_{\lambda}=x_{n_{\lambda}}$ and
$z_{\lambda}=z_{n_{\lambda}}$, as we wanted to show. 

Thus, if we take any $g\in
\mathrm{Lip}(X)$ with associated constant $k$, and $n$ as above, we have that 
\begin{displaymath}
\left|g^{\gamma X} \pl x_n \pr -g^{\gamma X} \pl z_n \pr \right| \leq k \ 
d \pl x_n , z_n \pr .
\end{displaymath}

Clearly this implies that $g^{\gamma X} \pl x_0 \pr = g^{\gamma X} \pl z_0 \pr $. By the definition of  $\gamma X$, we have $x_0 = z_0$, and we are done.

Therefore we conclude that $d(K_1, L_1)>0$.
This lets us  consider the function
\begin{displaymath}
f(x):=\mathrm{max}\left\{0,1-\frac{d(x,L_1)}{d(K_1, L_1)} \right\}
\end{displaymath}
for all $x\in X$, defined in a similar way as in 
Remark~\ref{macaya},
 which belongs
to $\mathrm{Lip}(X)$ and satisfies $0\leq f\leq 1$, $f(K_1)\equiv
0$, and $f(L_1)\equiv 1$. 
This proves the lemma.
\end{proof}

The next lemma is a Lipschitz version (with a similar proof) of the result given in
\cite[Lemma 3.4]{AF} in the context of uniformly continuous
functions.

\begin{lemma}\label{gedelta}
Let $X$ be a complete metric space and let $x\in \gamma X$. Then,
$x$ is a $G_{\delta}$-set in $\gamma X$ if and only if $x\in X$.
\end{lemma}

We close this section with a result concerning sums of Lipschitz functions that will be used in next sections.

\begin{lemma}\label{berano}
Let $(f_{n})$ be a sequence of functions in $\mathrm{Lip}(X,E)$ with
pairwise disjoint cozero sets and suppose that there exists a constant 
$M>0$ such that $L(f_{n})\leq M$ for all $n\in \mathbb{N}$. If
$f:=\sum_{n=1}^{\infty}f_{n}$ belongs to $C(X,E)$, then $f$ is a
Lipschitz function.
\end{lemma}

\begin{proof}
Let $x,y\in X$. Suppose first that $f(x)=f_{n_{0}}(x)$ and
$f(y)=f_{n_{0}}(y)$ for some $n_{0}\in \mathbb{N}$. Then
$\|f(x)-f(y)\|=\|f_{n_{0}}(x)-f_{n_{0}}(y)\|\leq M \ d(x,y)$.
Next assume that $f(x)=f_{n}(x)\neq 0$ and $f(y)=f_{m}(y)\neq 0$
with $n\neq m$. Then $\|f(x)-f(y)\|=\|f_{n}(x)-f_{m}(y)\|\leq
\|f_{n}(x)\|+\|f_{m}(y)\|=\|f_{n}(x)-f_{n}(y)\|+\|f_{m}(y)-f_{m}(x)\|\leq
2 M \ d(x,y)$. Consequently $L(f)\leq 2M$ and $f$ is a
Lipschitz function. 
\end{proof}

\section{Biseparating maps. Proofs}\label{sect4}

In this section we give the proofs of Theorems~\ref{imse2008} and \ref{segundano} and that of Proposition~\ref{kurda}, and
 some corollaries as well. We start with the notions of support point and support map.

\begin{definition}\label{jai-alai}
Let $T:\mathrm{Lip}(X,E)\rightarrow \mathrm{Lip}(Y,F)$ be a
biseparating map. A point $x\in  \gamma X$ is said to be a
\emph{support point} of $y\in Y$ if, for every neighborhood $U$ of
$x$ in $ \gamma X$, there exists $f\in \mathrm{Lip}(X,E)$ with
$\mathrm{coz}(f)\subset U$ such that $Tf(y)\neq 0$.
\end{definition}

\begin{remark}
For each $y\in Y$, the support point of $y\in Y$ exists and is unique  (see \cite[Lemma 4.3]{A3}). This fact
lets us define a map $h_{T}:  Y \rightarrow  \gamma X$ sending each
$y\in  Y$ to its support point $h_{T}(y)\in  \gamma X$. This map is 
usually called the \emph{support map} of $T$. If there is no chance of  confusion,  
we will denote it just by $h$ (instead of $h_T$).
\end{remark}

\begin{proposition}\label{hrmo}
Let $T:\mathrm{Lip}(X,E)\rightarrow \mathrm{Lip}(Y,F)$ be a
biseparating map. Then $h(Y) \subset X$ and 
$h :Y\rightarrow X$ is a homeomorphism.
\end{proposition}

\begin{proof}
In view of \cite[Lemma 4.7]{A3}, we can define the extension 
$\widetilde{h}:\gamma Y\rightarrow \gamma X$  of $h$. Besides,
taking into account Lemmas~\ref{trestres}, \ref{tres4}, and \ref{tresc}, we deduce that
$\widetilde{h}$ is a homeomorphism by applying \cite[Theorem 3.1]{A3}.
On the other hand, we have characterizated the points in $X$ as being the only
$G_{\delta}$-points in $\gamma X$ (see Lemma~\ref{gedelta}). Then, for each $y\in Y$,
$h (y)$ clearly belongs to $ X$ and $h:Y\rightarrow X$ is a 
a homeomorphism.
\end{proof}

\begin{lemma}\label{antto}
If $T:\mathrm{Lip}(X,E)\rightarrow \mathrm{Lip}(Y,F)$ is a
biseparating map and $f\in \mathrm{Lip}(X,E)$ satisfies $f\equiv 0$
on a neighborhood of $h (y)$, then $Tf\equiv 0$ on a neighborhood of $y$.
\end{lemma}

\begin{proof}
See \cite[Lemma 4.4]{A3}.
\end{proof}

\begin{lemma}\label{mazios}
Let $T:\mathrm{Lip}(X,E)\rightarrow \mathrm{Lip}(Y,F)$ be a
biseparating map. Let $f\in \mathrm{Lip}(X,E)$ and $y_{0}\in Y$ be such that
$f(h(y_{0}))=0$. Then $Tf(y_{0})=0$.
\end{lemma}

\begin{proof}
Let $(r_{n})$ be a sequence in $\mathbb{R}^{+}$ which converges to
$0$ and satisfies $2r_{n+1}<r_{n}$ for every $n\in \mathbb{N}$. We
set $B_{n}:=B(h(y_{0}),r_{n})$, $B_{n}^{2}:=B(h(y_{0}),2r_{n})$, and
$\varphi_{n}:=\varphi_{h(y_{0}),r_{n}}$ for each $n\in \mathbb{N}$,
where $\varphi_{h(y_{0}),r_{n}}$ is given as in Remark~\ref{macaya}.

\begin{claim}\label{tqhm}
Let $n, m \in \mathbb{N}$, $n \neq m$. Then 
$$\pl B_{2n}^{2} \backslash B_{2n+1} \pr \cap \pl B_{2m}^{2} \backslash B_{2m+1} \pr = \emptyset 
= \pl B_{2n-1}^{2} \backslash B_{2n} \pr \cap \pl B_{2m-1}^{2} \backslash B_{2m} \pr .$$
 
\end{claim}

The proof of Claim~\ref{tqhm} follows directly from the fact that, for all $k\in\mathbb{N}$,  $2r_{k+1}<r_{k}$, and
consequently $B_{k+1}^2 \subset B_{k}$.

\begin{claim}\label{pindar}
$L(f\varphi_{n})\leq 3 \ L(f)$ for all $n\in \mathbb{N}$.
\end{claim}

It is clear that $f\varphi_{n}\in \mathrm{Lip}(X,E)$ for all $n\in \mathbb{N}$.
Now, by definition of $\varphi_{n}$,
$\mathrm{coz}(f\varphi_{n})\subset B_{n}^{2}$, and if $x\in B_{n}^{2}$, then
$\|f(x)\|=\|f(x)-f(h(y_{0}))\|\leq L(f) \ d(x,h(y_{0}))< 2r_{n} \ L(f)$. 
Consequently, if $x,y\in B_{n}^{2}$,
\begin{eqnarray*}
\|(f\varphi_{n})(x)-(f\varphi_{n})(y)\|&\leq&
\|f(x)\|\left|\varphi_{n}(x)-\varphi_{n}(y)\right|+\left|\varphi_{n}(y)\right|
\|f(x)-f(y)\|\\
&\leq& 2r_{n} \ L(f) \ (1/r_{n}) \ d(x,y)+L(f) \ d(x,y)\\
&=& 3 \ L(f) \ d(x,y).
\end{eqnarray*}
Besides, if $x\in B_{n}^{2}$ and $y\notin B_{n}^{2}$,
\begin{eqnarray*}
\|(f\varphi_{n})(x)-(f\varphi_{n})(y)\|\leq 2r_{n} \ L(f) \ (1/r_{n}) \ d(x,y)=2 \ L(f) \ d(x,y). 
\end{eqnarray*}
Thus Claim~\ref{pindar} is proved.

\smallskip
Next we consider the function $g:=f\varphi_{1}$, and
define  $g_{1}:=\sum_{n=1}^{\infty}f \ (\varphi_{2n}-\varphi_{2n+1})$
and $g_{2}:=\sum_{n=1}^{\infty}f \ (\varphi_{2n-1}-\varphi_{2n})$.
It is obvious that
$g=g_{1}+g_{2}$, and since $f(h(y_{0}))=0$, we see that $g_{1}(h(y_{0}))=0$
and $g_{2}(h(y_{0}))=0$. This implies that both $g_1$ and $g_2$ are continuous. 
Taking into account Claim~\ref{pindar}, $L(f(\varphi_{n}-\varphi_{n+1}))\leq
L(f\varphi_{n})+L(f\varphi_{n+1})\leq 6 \ L(f)$ for all $n\in
\mathbb{N}$. Besides, since
$\mathrm{coz}(\varphi_{2n}-\varphi_{2n+1})\subset
B_{2n}^{2}\backslash B_{2n+1}$, we deduce from Claim~\ref{tqhm} that
\begin{displaymath}
\mathrm{coz}(\varphi_{2n}-\varphi_{2n+1})\cap
\mathrm{coz}(\varphi_{2m}-\varphi_{2m+1})=\emptyset
\end{displaymath}
whenever $n\neq m$. Applying Lemma~\ref{berano}, we conclude that $g_{1}$ (and
similarly $g_{2}$) belongs to $\mathrm{Lip}(X,E)$. Besides, $g\equiv f$ on $B_{1}$, and
by Lemma~\ref{antto}, $Tg(y_{0})=Tf(y_{0})$. Therefore, to see that $Tf(y_{0})=0$,
it is enough to prove that $Tg_{1}(y_{0})=0$ and $Tg_{2}(y_{0})=0$.

\begin{claim}\label{nasio}
Given $n_0 \in \mathbb{N}$, $$\mathrm{cl}_{X}\pl \mathrm{coz}\pl g_{1} \pr \pr \subset \mathrm{cl}_{X}\pl B_{2n_0}^2 \pr \cup  \bigcup_{n=1}^{n_0-1} \mathrm{cl}_{X}\pl  B_{2n}^2  \setminus B_{2n+1} \pr   .$$
\end{claim}

To see this, notice  that  
\begin{eqnarray*}
\mathrm{coz}\left(\sum_{n=n_{0}}^{\infty} \varphi_{2n}-\varphi_{2n+1} \right) &\subset&
\bigcup_{n=n_{0}}^{\infty}\mathrm{coz}(\varphi_{2n}-\varphi_{2n+1}) \\
&\subset& B_{2n_{0}}^{2},
\end{eqnarray*}
and that $\mathrm{coz} \pl \varphi_{2n}-\varphi_{2n+1} \pr \subset   B_{2n}^{2} \setminus B_{2n +1}
$ for $n < n_0$.

\medskip

If we consider, for each $n\in \mathbb{N}$, a point $y_{n}\in
h^{-1}(B_{2n-1})\backslash \mathrm{cl}_{Y}h^{-1}(B_{2n}^{2})$, then the sequence 
$(y_{n})$ converges to $y_{0}$ because $\cap_{n=1}^{\infty}B_{n}=\{h(y_{0})\}$ 
and $h$ is a homeomorphism.

\begin{claim}\label{inbi}
$h(y_{n})\notin \mathrm{cl}_{X}\pl \mathrm{coz}(g_{1}) \pr$ for all $n\in \mathbb{N}$.
\end{claim}

Let us prove the claim. Fix $n_{0}\in \mathbb{N}$. 
It is clear by construction that
$h(y_{n_{0}})\notin \mathrm{cl}_{X}(B_{2n_{0}}^{2})$ and that,
if $n<n_{0}$, then 
 $h(y_{n_{0}})\in B_{2n_{0}-1}\subseteq B_{2n+1}$, that is, $h(y_{n_{0}})\notin \mathrm{cl}_{X} \pl B_{2n}^2 \setminus  B_{2n+1} \pr$. Therefore Claim~\ref{inbi} follows from  Claim~\ref{nasio}. 
 
\medskip

Finally, since $h(y_{n})\notin \mathrm{cl}_{X}(\mathrm{coz}(g_{1}))$ for all $n\in \mathbb{N}$, 
then $g_1 \equiv 0$ on a neighborhood of $h(y_n)$.
Applying Lemma~\ref{antto}, 
$Tg_{1}(y_{n})=0$ for all $n\in \mathbb{N}$ , and by continuity, we conclude that 
$Tg_{1}(y_{0})=0$. In the same way it can be proved that $Tg_{2}(y_{0})=0$.
\end{proof}

\begin{proposition}\label{sardinia}
Let $T:\mathrm{Lip}(X,E)\rightarrow \mathrm{Lip}(Y,F)$ be a
biseparating map. For each $y \in Y$, there exists a linear and bijective map $J y  : E \ra  F$ such that
$$Tf(y)=(J y)(f(h(y)))$$
for all $f\in \mathrm{Lip}(X,E)$ and $y\in Y$.
\end{proposition}

\begin{proof}
For $y\in Y$ and $f\in \mathrm{Lip}(X,E)$ fixed, consider the function $g:=f-\widehat{f(h(y))}\in
\mathrm{Lip}(X,E)$. Clearly
$g(h(y))=0$, and by Lemma ~\ref{mazios}, $Tg(y)=0$. Consequently 
$Tf(y)=T\widehat{f(h(y))}(y)$ for all $f\in \mathrm{Lip}(X,E)$ and
$y\in Y$. Next, we define $J y :E\rightarrow F$ as $(J y)(\mathbf{e}):=T\widehat{\mathbf{e}}(y)$ 
for all $\mathbf{e}\in E$, which is linear and bijective (see \cite[Theorem 3.5]{A2}). We easily see that 
$T$ has the desired representation.
\end{proof}

\begin{remark}\label{sinbox}
Notice that, if $T:\mathrm{Lip}(X,E)\rightarrow \mathrm{Lip}(Y,F)$ is a
biseparating map, $T^{-1}:\mathrm{Lip}(Y,F)\rightarrow \mathrm{Lip}(X,E)$
is also biseparating, so there exist a homeomorphism $h_{T^{-1}}:X\rightarrow
Y$ and a map $K x :F\rightarrow E$ for all $x\in X$ such that
\begin{center}
$T^{-1}g(x)=(K x )(g(h_{T^{-1}}(x)))$ 
\end{center}
for all $g\in \mathrm{Lip}(Y,F)$ and $x\in X$.
Besides, it is not difficult to check that $h_{T^{-1}}\equiv h_{T}^{-1}$ 
(see Claim 1 in the proof of the Theorem 3.1 in \cite{A3}).
\end{remark}

\begin{lemma}\label{acp}
Let $T:\mathrm{Lip}(X,E)\rightarrow \mathrm{Lip}(Y,F)$ be a
biseparating map. Then $\mathrm{inf}\{\| (J y)(\mathbf{e})\|:y\in Y\}>0$
for each non-zero $\mathbf{e}\in E$.
\end{lemma}

\begin{proof}
Suppose this is not true. Then there exist $(y_{n})$ in $Y$ and $\mathbf{e}\in E$ with
$\|\mathbf{e}\|=1$ such that
$\| \pl J y_{n} \pr (\mathbf{e})\|<1/n^{3}$ for each $n\in \mathbb{N}$.

If we assume first that there exists a limit point
$y_{0}\in Y$ of $\{y_{n}:n\in \mathbb{N}\}$, then we can consider a
subsequence $(y_{n_{k}})$ of $(y_{n})$ converging to $y_{0}$, so that
$\| \pl J y_{0} \pr (\mathbf{e})\|=0$, which is absurd since $J y_{0}$ is inyective.

Therefore, there exists $r>0$ such that $d(y_{n},y_{m})>r$
whenever $n\neq m$. Also, on the one hand,
$\ql T^{-1}(T\widehat{\mathbf{e}}) \qr (h(y_{n}))=\widehat{\mathbf{e}}(h(y_{n}))=\mathbf{e}$ for all
$n\in \mathbb{N}$, and on the other hand, by Remark~\ref{sinbox},
$\ql T^{-1}(T\widehat{\mathbf{e}}) \qr (h(y_{n}))= (K h(y_{n}) ) (T\widehat{\mathbf{e}}(y_{n}))$.
Consequently $\| \pl K h(y_{n}) \pr (T\widehat{\mathbf{e}}(y_{n}))\|=\vc \mathbf{e} \vd =1$ for each $n\in \mathbb{N}$.
If we take $\mathbf{f}_{n}\in F$ defined as
$\mathbf{f}_{n}:=T\widehat{\mathbf{e}}(y_{n})/\|T\widehat{\mathbf{e}}(y_{n})\|$ for each $n\in \mathbb{N}$,
it is clear that $\|\mathbf{f}_{n}\|=1$ and
$$\| \pl K h(y_{n}) \pr (\mathbf{f}_{n})\|=
(1/\|T\widehat{\mathbf{e}}(y_{n})\|)\| \pl K h(y_{n}) \pr (T\widehat{\mathbf{e}}(y_{n}))\|>n^{3}.$$

Next, we define, in a similar way as in 
Remark~\ref{macaya},
\begin{displaymath}
\psi_{y_{n},r/3}(y):=\mathrm{max}\left\{0,1- \frac{3 \ d(y,y_{n})}{r} \right\}
\end{displaymath}
for all $y\in Y$ and $n\in \mathbb{N}$ (denoted for short $\psi_{n}$) which
belongs to $\mathrm{Lip}(Y)$, and finally, we consider the function
\begin{displaymath}
g:=\sum_{n=1}^{\infty}\frac{\psi_{n}\mathbf{f}_{n}}{n^{2}}.
\end{displaymath}
It is immediate to see that $\|\psi_{n}\mathbf{f}_{n}/n^{2}\|_{\infty}\leq 1/n^{2}$ and
$L(\psi_{n}\mathbf{f}_{n}/n^{2})=(\|\mathbf{f}_{n}\|/n^{2})L(\psi_{n})=3/\pl rn^{2} \pr$  for 
all $n\in \mathbb{N}$, which lets us conclude by Lemma~\ref{berano} that 
 $g$ belongs
to $\mathrm{Lip}(Y,F)$.

It is apparent that $g(y_{n})=\mathbf{f}_{n}/n^{2}$, and applying 
Lemma~\ref{mazios} for the biseparating map $T^{-1}$, we deduce that 
$T^{-1}g(h(y_{n}))=(1/n^{2})T^{-1}\widehat{\mathbf{f}_{n}}(h(y_{n}))$.
Consequently, $\vc T^{-1}g(h(y_{n})) \vd =
(1/n^{2})\vc \pl K h(y_{n}) \pr (\mathbf{f}_{n})\vd >n$ for all $n\in
\mathbb{N}$, which contradicts the fact that $T^{-1}g$ is bounded.
\end{proof}

\begin{proof}[Proof of Proposition~\ref{kurda}]
Suppose on the contrary that there exist sequences $(y_n)$ in $Y $ and $ \pl \mathbf{e}_n \pr$ in $E$
with $\vc \mathbf{e}_n \vd=1$ and $\vc T \widehat{\mathbf{e}_n} (y_n) \vd > n^2$ for every $n \in \mathbb{N}$. Take $\mathbf{f} \in F$ with $\vc \mathbf{f} \vd =1$. By Lemma~\ref{acp} there exists $M >0$ such that $\vc  T^{-1} \widehat{\mathbf{f}} (h(y_n) ) \vd > M$ for every $n$. Consider a sequence $(r_n)$ in $(0,1)$ such that $B(y_n, r_n) \cap B(y_m, r_m) = \emptyset$ whenever $n \neq m$ (this can be done by taking a subsequence of $(y_n)$ if necessary).
Without loss of generality we may also assume that $(r_n)$ is decreasing and converging to $0$.

We define, for each $n\in \mathbb{N}$,
\begin{displaymath}
\xi_{n}(y):=\mathrm{max}\{0,r_{n}-d(y,y_{n})\}
\end{displaymath}
for all $y \in Y$, which belongs to $\mathrm{Lip}(Y)$ and satisfies
$\xi_{n}(y_{n}) =r_{n}$,
$\mathrm{coz}(\xi_{n})=B(y_{n},r_{n})$,
$\|\xi_{n}\|_{\infty}=r_{n}$, and $L(\xi_{n})=1$. Finally, we consider the function
\begin{displaymath}
g:=\sum_{n=1}^{\infty} \xi_{n} \mathbf{f}.
\end{displaymath}
The fact  that $g$ belongs to $\mathrm{Lip}(Y,F)$ follows from  Lemma~\ref{berano}. Now let $f:= T^{-1} g $. 
It is clear that $f: = \sum_{n=1}^{\infty} T^{-1} \pl \xi_{n} \mathbf{f} \pr$. Consequently, if for each $n \in \mathbb{N}$, we define $f_n (x) := \vc T^{-1} \pl \xi_{n} \mathbf{f} \pr (x) \vd$ ($x \in X$), then  $f_0 := \vc f \vd = \sum _{n=1}^{\infty}
f_n$ belongs to $\mathrm{Lip}(X)$ and $f_0 (h(y_n)) \ge M r_n$ for every $n \in \mathbb{N}$. Therefore $f_0' :=  \sum _{n=1}^{\infty}
f_n \mathbf{e_n}$ belongs to $\mathrm{Lip}(X,E)$. Finally $\vc Tf_0' (y_n) \vd \ge Mr_n n^2$, and it is easily seen that
$L \pl T f_n \mathbf{e_n} \pr \ge M n^2$, for every $n \in \mathbb{N}$. We conclude that $Tf_0'$ does not belong to 
$\mathrm{Lip}(Y,F)$, which is absurd.

Now, the fact that each $y \in Y_d$ is isolated follows easily.
\end{proof}

\begin{remark}
We will  use later the fact that, since $Y_{d}$ is a finite set of isolated points and $h$ is a 
homeomorphism, then  $d(X\setminus h(Y_{d}), h(Y_{d}))>0$.
\end{remark}

The restriction to $X\setminus h(Y_{d})$ (respectively, $Y\setminus Y_{d}$) of a function  
$f \in \mathrm{Lip}(X,E)$ (respectively, $ f \in \mathrm{Lip}(Y,F)$), is obviously   a bounded Lipschitz function, which 
will be denoted by $f_d$. The converse is also true, that is, we can obtain a Lipschitz function as an extension of 
an element of $\mathrm{Lip}(X\setminus h(Y_{d}),E)$, as it is done in the next lemma.

\begin{lemma}\label{lute}
Let $f \in \mathrm{Lip}(X\setminus h(Y_{d}),E)$. Then  
the function
\begin{displaymath}
f^d (x):= \left\{ \begin{array}{ll}
      f (x) & \mbox{{\rm if} } x\in X\setminus h(Y_{d})\\
      0 & \mbox{{\rm if} } x\in h(Y_{d}) 
      \end{array} \right.
\end{displaymath}
belongs to $\mathrm{Lip}(X,E)$.
\end{lemma}

\begin{proof}
Since $h(Y_{d})$ is a finite set of isolated points, 
 $f^d $ is clearly a continuous function. Besides, if we consider $x_{1}\in X\setminus h(Y_{d})$ 
and $x_{2}\in h(Y_{d})$,
\begin{displaymath}
\frac{\|f^d \pl x_{1} \pr - f^d (x_{2})\|}{d(x_{1},x_{2})}\leq \frac{\|f (x_{1})\|}{ d(X\setminus h(Y_{d}), h(Y_{d} ))}\leq
\frac{\|f\|_{\infty}}{d(X\setminus h(Y_{d}),h(Y_{d}))}.
\end{displaymath}
Therefore $$L \pl f^d \pr \leq \mathrm{max}\tl L(f),\frac{\|f\|_{\infty}}{ d(X\setminus h(Y_{d}), h(Y_{d}))}\tr <\infty,$$
which implies that $f^d \in \mathrm{Lip}(X,E)$.
\end{proof}

\begin{proof}[Proof of Theorem~\ref{segundano}]
By definition of $T_d$ and Lemma~\ref{lute} (see also the comment before it), we clearly see that $$T_d \pl f \pr = \pl T f^d \pr_d$$ for all $f\in \mathrm{Lip}(X\setminus h(Y_{d}),E)$, so $T_d$ is well defined and it is biseparating. To prove 
that $T_d$ is continuous, we will see that 
given a sequence $(f_{n})$ in $\mathrm{Lip}(X\setminus h(Y_{d}),E)$ 
converging to $0$ and such that $\pl T_d f_{n} \pr$ converges to 
$g \in \mathrm{Lip}(Y\setminus Y_{d},F)$, we have $g \equiv 0$.

If we consider, for each $n\in \mathbb{N}$, the extension
$f_{n}^d$ of $f_{n}$ given in Lemma~\ref{lute}, we can show that 
\begin{eqnarray*}
\vc f_{n}^d \vd_{L}
&\leq& \max \left\{\vc f_{n} \vd_{\infty},
\max \tl L(f_{n}), \frac{\vc f_{n} \vd_{\infty}}{d(X\setminus h(Y_{d}),h(Y_{d}))} \tr\right\} \\
&\leq& \vc f_{n} \vd_{L} \max \tl 1,1/d(X\setminus h(Y_{d}),h(Y_{d}))\tr ,
\end{eqnarray*}
which allows us to deduce that $(f_{n}^d)$ converges to $0$.
By continuity, if we fix $y\in Y\setminus Y_{d}$, the sequence
$\pl \pl J y \pr  \pl  f_{n}^d (h(y)) \pr \pr$
converges to $0$. Besides, since $Tf_{n}^d (y)=T_d f_{n}(y)$, 
we conclude that 
$\pl T_d f_{n}(y) \pr$
converges to $0$.

On the other hand, $\vc T_d f_{n} (y)- g(y) \vd
\leq \vc T_d f_{n}- g \vd_{L}$ for each $n\in \mathbb{N}$, and 
as $(T_d f_{n})$ converges to $g$, we deduce that 
$\pl T_d f_{n} (y) \pr$ converges to $g(y)$.
Combined with the above, $g(y)=0$ for all $y\in Y\backslash Y_{d}$.
\end{proof}

The proof of the two following results is now immediate.

\begin{corollary}\label{linda1}
Suppose that $E$ and $F$ are complete and let $T:\mathrm{Lip}(X,E)\rightarrow \mathrm{Lip}(Y,F)$ be a
biseparating map. If $Y$ has no isolated points, then $T$ is continuous.
\end{corollary}

\begin{corollary}\label{linda}
Let $T:\mathrm{Lip}(X,E)\rightarrow \mathrm{Lip}(Y,F)$ be a
biseparating map. If $E$ has finite dimension, then $F$ has the same
dimension as $E$ and $T$ is continuous.
\end{corollary}

\begin{proposition}\label{pamela}
 Let $T:\mathrm{Lip}(X,E)\rightarrow \mathrm{Lip}(Y,F)$ be a
biseparating map. Then $h:Y\rightarrow X$ is a bi-Lipschitz map.
\end{proposition}

\begin{proof}
Associated to $T$, we define a linear map $S:\mathrm{Lip}(X)\rightarrow \mathrm{Lip}(Y)$.
For $f \in \mathrm{Lip}(X)$, define $$Sf (y) := f(h(y))$$ for every
$y \in Y$. It is obvious that $Sf$ is a continuous bounded function on $Y$. Next we are going to see that
it is also Lipschitz. It is clear that it is enough to prove it in the case when $f \ge 0$.

Fix any $\mathbf{e} \neq 0$ in $E$. By Lemma~\ref{acp}, we know that there exists $M >0$ such that $\vc T \widehat{\mathbf e} (y) \vd \ge M$  for every $y \in Y$, so the map $y \mapsto 1/\vc T \widehat{\mathbf e} (y) \vd $ belongs to $\mathrm{Lip}(Y)$. On the other hand, taking into account that $f \ge 0$, we have that  for $y , y' \in Y$
\begin{eqnarray*}
\va Sf (y)  \vc T \widehat{\mathbf{e}} (y) \vd  - Sf (y')  \vc T \widehat{\mathbf{e}}  (y') \vd \vb
 &=& \va \vc \pl J y \pr (f (h(y)) \mathbf{e} ) \vd - \vc \pl J y' \pr (f (h(y')) \mathbf{e} ) \vd  \vb \\
&\le&   \vc (J y) (f (h(y)) \mathbf{e} )  -  (J y') (f (h(y')) \mathbf{e} ) \vd  \\
&=&  \vc T (f \mathbf{e}) (y) - T (f \mathbf{e}) (y') \vd \\
&\le& L (T (f \mathbf{e})) \ d(y,y').
\end{eqnarray*}
We deduce that $Sf$ is Lipschitz. A similar process can be done with the map $T^{-1}$, and we
conclude that $S:\mathrm{Lip}(X)\rightarrow \mathrm{Lip}(Y)$ is bijective and biseparating.

Next we prove that $h$ is Lipschitz. 
Let $K_0 := \max \tl 1, \mathrm{diam}(X) \tr$. 
We  take $y,y' \in Y$ and define
$f_{1}(x):=d(h(y),x)$ for all $x\in X$. Clearly $f_{1}$ belongs to $\mathrm{Lip}(X)$ and,
since  $S$ is continuous (see Corollary~\ref{linda}), it is not difficult to see that
\begin{displaymath}
\frac{\vc S f_{1}(y)- S f_{1}(y')\vd}{d(y,y')}
\leq \vc S f_{1} \vd_{L}\leq
\vc S \vd \vc f_{1} \vd_{L} \leq K_0 \vc S \vd.
\end{displaymath}
On the other hand,
$Sf_{1}(y) =0$ and $S f_{1}(y') =  d(h(y),h(y'))$.
Then, replacing in the above inequality,
$$d(h(y),h(y'))\leq 
K_0 \vc S \vd \ d(y,y'),$$
and we are done.

Moreover, $h^{-1}$ is also Lipschitz because $h^{-1} = h_{T^{-1}}$ (see Remark~\ref{sinbox})
\end{proof}

\begin{proof}[Proof of Theorem~\ref{imse2008}]
It follows immediately from Propositions~\ref{hrmo}, \ref{sardinia} and \ref{pamela}.
\end{proof}

Taking into account Theorem~\ref{imse2008}, Lemma~\ref{acp}, and Corollary~\ref{linda}, we can give the general form of biseparating maps in 
the scalar-valued case (see also Theorem~\ref{gerria} and Corollary~\ref{labaina}). Of course it also applies to algebra isomorphisms.

\begin{corollary}\label{tepeyac}
Let $T:\mathrm{Lip}(X)\rightarrow \mathrm{Lip}(Y)$ be a biseparating map. Then $T$ is continuous and there exist a bi-Lipschitz
homeomorphism $h:Y\rightarrow X$ and a nonvanishing function
$\tau\in \mathrm{Lip}(Y)$ such that
\begin{center}
$Tf(y)=\tau(y)f(h(y))$
\end{center}
for every $f\in \mathrm{Lip}(X)$ and $y\in Y$.
\end{corollary}

\begin{corollary}\label{nadena}
Let $I:\mathrm{Lip}(X)\rightarrow \mathrm{Lip}(Y)$ be an algebra isomorphism. Then $I$ is continuous and there exists a bi-Lipschitz
homeomorphism $h:Y\rightarrow X$such that
\begin{center}
$Tf(y)= f(h(y))$
\end{center}
for every $f\in \mathrm{Lip}(X)$ and $y\in Y$. 
\end{corollary}

\section{Separating maps. Proof of Theorem~\ref{gerria}}\label{sect5}

In this section we give the proof of Theorem~\ref{gerria} and the representation of    bijective separating maps in the scalar setting when $Y$ is compact.

\begin{proof}[Proof of Theorem~\ref{gerria}]
Let $f,g\in \mathrm{Lip}(X)$ be such that $\mathrm{coz}(f)\cap
\mathrm{coz}(g)\neq \emptyset$, that is, there exists $x_{0}\in X$
satisfying $f(x_{0})\neq 0$ and $g(x_{0})\neq 0$. Since $T$ is onto,
$Tk\equiv  1$ 
for some $k\in \mathrm{Lip}(X)$, and we can take
$\alpha,\beta\in \mathbb{K}$ such that $(\alpha f+k)(x_{0})=0$ and 
$(\beta g+k)(x_{0})=0$. We denote $l:=\alpha f+k$.

Let $(r_{n})$, $B_{n}$, $B_{n}^{2}$,  and $\varphi_{n}$ be as in the proof of
Lemma~\ref{mazios} (where $h(y_{0})$ is replaced by $x_{0}$); indeed, we closely follow that proof. 
Now, we take $y_{n}\in \mathrm{coz}(T(\varphi_{n}-\varphi_{n+1}))$
for each $n\in \mathbb{N}$. By the compactness of $Y$, $\{y_{n}:n\in
\mathbb{N}\}$ has a limit point $y_{0}$ in $Y$. Then, we can consider a
subsequence $(y_{n_{i}})$ of $(y_{n})$ converging to $y_{0}$ whose
indexes satisfy $\va n_{i}-n_{j} \vb \geq 3$ whenever $i\neq j$.

We claim that $Tl(y_{0})=0$. To prove it, we define
\begin{displaymath}
l_{1}:=\sum_{k=1}^{\infty}l(\varphi_{n_{2k}-1}-\varphi_{n_{2k}+2})
\end{displaymath}
and $l_{2}:=l -l_{1}$, and we will see that $Tl_{1}(y_{0})=0$ and
$Tl_{2}(y_{0})=0$ (in the rest of the proof we will set 
$\xi_{k}:=\varphi_{n_{2k}-1}-\varphi_{n_{2k}+2}$ for every $k\in \mathbb{N}$).

First, we will check that $l_{1}$ and $l_{2}$ are both Lipschitz
functions. As in  Claim~\ref{pindar} in the proof of  Lemma~\ref{mazios}, we know that
$L(l\varphi_{n})\leq 3L(l)$ for all $n\in \mathbb{N}$. Consequently
$L(l\xi_{k})\leq L(l\varphi_{n_{2k}-1})+L(l\varphi_{n_{2k}+2})\leq 6L(l)$ for all
$k\in \mathbb{N}$. Since
$\mathrm{coz}(\xi_{k})\cap
\mathrm{coz}(\xi_{j})=\emptyset$ if
$k\neq j$, by Lemma~\ref{berano} we conclude that $l_{1}\in \mathrm{Lip}(X)$,
and then $l_{2}$ also belongs to $\mathrm{Lip}(X)$.

Now, we will see that $Tl_{1}(y_{n_{2k-1}})=0$ for all $k\in
\mathbb{N}$. Fix $k_{0}\in\mathbb{N}$ and consider
$y_{n_{2k_{0}-1}}$. It is not difficult to see that
$\mathrm{coz}(\varphi_{n_{2k_{0}-1}}-\varphi_{n_{2k_{0}-1}+1}) \subset B_{n_{2k_{0}-1}}^{2}\backslash B_{n_{2k_{0}-1}+1}$ and that, for every $k\in
\mathbb{N}$,
$\mathrm{coz}(\xi_k)\subset
B_{n_{2k}-1}^{2}\backslash B_{n_{2k}+2}$, so 
$$\mathrm{coz}(\varphi_{n_{2k_{0}-1}}-\varphi_{n_{2k_{0}-1}+1})\cap
\mathrm{coz}(\xi_{k})=\emptyset,$$
which allows us to deduce that
$$\mathrm{coz}\left(\varphi_{n_{2k_{0}-1}}-\varphi_{n_{2k_{0}-1}+1}\right)\bigcap
\mathrm{coz}\left(\sum_{k=1}^{\infty} l\xi_{k}\right)=\emptyset.$$
Next, since $T$ is a separating map,
$$\mathrm{coz}\left(T(\varphi_{n_{2k_{0}-1}}-\varphi_{n_{2k_{0}-1}+1})\right)\bigcap
\mathrm{coz}\left(T\left(\sum_{k=1}^{\infty} l\xi_{k}\right)\right)=\emptyset,$$
and we conclude that
$Tl_{1}(y_{n_{2k_{0}-1}})=0$ because $y_{n_{2k_{0}-1}}\in
\mathrm{coz}(T(\varphi_{n_{2k_{0}-1}}-\varphi_{n_{2k_{0}-1}+1})) $.
By continuity, it is clear that $Tl_{1}(y_{0})=0$.

On the other hand, if $x\in \mathrm{coz}(\varphi_{n_{2k}}-\varphi_{n_{2k}+1})=
B_{n_{2k}}^{2}\backslash B_{n_{2k}+1}\subset B_{n_{2k}-1}\backslash B_{n_{2k}+2}^{2}$, 
then $\xi_{k}(x)=1$. This fact allows us to deduce that
$\mathrm{coz}(\varphi_{n_{2k}}-\varphi_{n_{2k}+1})\cap
\mathrm{coz}(l_{2})=\emptyset$, and consequently
$\mathrm{coz}(T(\varphi_{n_{2k}}-\varphi_{n_{2k}+1}))\cap
\mathrm{coz}(Tl_{2})=\emptyset$. For this reason
$Tl_{2}(y_{n_{2k}})=0$ for all $k\in \mathbb{N}$, and as above we conclude that
$Tl_{2}(y_{0})=0$.

Therefore
$0=Tl(y_{0})=T(\alpha f+k)(y_{0})=\alpha Tf(y_{0})+1$, which
implies that $Tf(y_{0})\neq 0$. The same reasoning can be applied to 
the function $\beta g+k$ and we obtain that $Tg(y_{0})\neq 0$. Then, 
we deduce that $\mathrm{coz}(Tf)\cap \mathrm{coz}(Tg)\neq \emptyset$, 
and $T^{-1}$ is separating.

The fact that $T$ is continuous follows from Corollary~\ref{linda}.
\end{proof}

\begin{corollary}\label{labaina}
Let $T:\mathrm{Lip}(X)\rightarrow \mathrm{Lip}(Y)$ be a bijective
and separating map. If $Y$ is compact, then there exist a bi-Lipschitz
homeomorphism $h:Y\rightarrow X$ and a nonvanishing function
$\tau\in \mathrm{Lip}(Y)$ such that
\begin{center}
$Tf(y)=\tau(y)f(h(y))$
\end{center}
for every $f\in \mathrm{Lip}(X)$ and $y\in Y$.
\end{corollary}

\begin{proof}
Immediate by Theorem~\ref{gerria} and Corollary~\ref{tepeyac}.
\end{proof}

\end{document}